\documentclass[]{amsart}
\usepackage{amsmath, amsthm, amscd, amsfonts, amssymb, graphicx, color}
\usepackage[bookmarksnumbered, colorlinks, plainpages]{hyperref}

\theoremstyle{definition}

\theoremstyle{remark}

\numberwithin{equation}{section}

\begin{document}
\setcounter{page}{1}
\begin{center}
{\bf  ARENS REGULARITY  AND FACTORIZATION PROPERTY}
\end{center}

\title[]{}

\author[]{KAZEM HAGHNEJAD AZAR   }

\address{}

\dedicatory{}

\subjclass[2000]{46L06; 46L07; 46L10; 47L25}

\keywords {Arens regularity, bilinear mapping,  topological
center, module action, factorization,  weakly compact}

\begin{abstract}
In this paper,  we  study the Arens regularity properties of module actions and we extend some proposition from Baker, Dales, Lau and others into general situations.  For Banach $A-bimodule$ $B$, let $Z_1(A^{**})$,  ${Z}^\ell_{B^{**}}(A^{**})$ and  ${Z}^\ell_{A^{**}}(B^{**})$ be the topological centers of second dual of Banach algebra $A$,  left module action $\pi_\ell:~A\times B\rightarrow B$ and  right module action $\pi_r:~B\times A\rightarrow B$, respectively. We establish some relationships between them and factorization properties of $A^*$ and $B^*$. We  search some necessary and sufficient conditions for factorization of $A^*$, $B$ and $B^*$ with some results in group algebras. We extend the definitions of the left and right multiplier for module actions.

\end{abstract}\maketitle

\section{\bf  Preliminaries and
Introduction }
\noindent In  1951 Arens  shows that  the second dual $A^{**}$ of Banach algebra $A$ endowed with the either Arens multiplications is a Banach algebra, see [1]. The constructions of the two Arens multiplications in $A^{**}$ lead us to definition of topological centers for $A^{**}$ with respect to both Arens multiplications. The topological centers of Banach algebras, module actions and applications of them  were introduced and discussed in [3, 5, 6, 9, 15, 16, 17, 18, 19, 24, 25]. In this paper, we  extend some problems from [3, 5, 6, 16,  22] to the general criterion on module actions with some applications in  group algebras. \\
Now we introduce some notations and definitions that we used
in this paper.\\
\noindent Throughout  this paper, $A$ is  a Banach algebra and $A^*$,
$A^{**}$, respectively, are the first and second dual of $A$.
 For $a\in A$
 and $a^\prime\in A^*$, we denote by $a^\prime a$
 and $a a^\prime$ respectively, the functionals on $A^*$ defined by $\langle a^\prime a,b\rangle= \langle  a^\prime,ab\rangle=a^\prime(ab)$ and $ \langle  a a^\prime,b\rangle= \langle  a^\prime,ba\rangle=a^\prime(ba)$ for all $b\in A$.
 The Banach algebra $A$ is embedded in its second dual via the identification
 $ \langle  a,a^\prime\rangle$ - $ \langle  a^\prime,a\rangle$ for every $a\in
A$ and $a^\prime\in
A^*$.
  We denote the set   $\{a^\prime a:~a\in A~ and ~a^\prime\in
  A^*\}$ and
  $\{a a^\prime:~a\in A ~and ~a^\prime\in A^*\}$ by $A^*A$ and $AA^*$, respectively, clearly these two sets are subsets of $A^*$.\\
 The extension of bilinear maps on normed space and the concept of regularity of bilinear maps were studied by [1, 2, 5, 6, 9]. We start by recalling these definitions as follows.\\
 Let $X,Y,Z$ be normed spaces and $m:X\times Y\rightarrow Z$ be a bounded bilinear mapping. Arens in [1] offers two natural extensions $m^{***}$ and $m^{t***t}$ of $m$ from $X^{**}\times Y^{**}$ into $Z^{**}$ as following\\
1. $m^*:Z^*\times X\rightarrow Y^*$,~~~~~given by~~~$ \langle  m^*(z^\prime,x),y\rangle= \langle  z^\prime, m(x,y)\rangle$ ~where $x\in X$, $y\in Y$, $z^\prime\in Z^*$,\\
 2. $m^{**}:Y^{**}\times Z^{*}\rightarrow X^*$,~~given by $ \langle  m^{**}(y^{\prime\prime},z^\prime),x\rangle= \langle  y^{\prime\prime},m^*(z^\prime,x)\rangle$ ~where $x\in X$, $y^{\prime\prime}\in Y^{**}$, $z^\prime\in Z^*$,\\
3. $m^{***}:X^{**}\times Y^{**}\rightarrow Z^{**}$,~ given by~ ~ ~$ \langle  m^{***}(x^{\prime\prime},y^{\prime\prime}),z^\prime\rangle$ $= \langle  x^{\prime\prime},m^{**}(y^{\prime\prime},z^\prime)\rangle$ ~\\where ~$x^{\prime\prime}\in X^{**}$, $y^{\prime\prime}\in Y^{**}$, $z^\prime\in Z^*$.\\
The mapping $m^{***}$ is the unique extension of $m$ such that $x^{\prime\prime}\rightarrow m^{***}(x^{\prime\prime},y^{\prime\prime})$ from $X^{**}$ into $Z^{**}$ is $weak^*-to-weak^*$ continuous for every $y^{\prime\prime}\in Y^{**}$, but the mapping $y^{\prime\prime}\rightarrow m^{***}(x^{\prime\prime},y^{\prime\prime})$ is not in general $weak^*-to-weak^*$ continuous from $Y^{**}$ into $Z^{**}$ unless $x^{\prime\prime}\in X$. Hence the first topological center of $m$ may  be defined as following
$$Z_1(m)=\{x^{\prime\prime}\in X^{**}:~~y^{\prime\prime}\rightarrow m^{***}(x^{\prime\prime},y^{\prime\prime})~~is~~weak^*-to-weak^*~~continuous\}.$$
Let now $m^t:Y\times X\rightarrow Z$ be the transpose of $m$ defined by $m^t(y,x)=m(x,y)$ for every $x\in X$ and $y\in Y$. Then $m^t$ is a continuous bilinear map from $Y\times X$ to $Z$, and so it may be extended as above to $m^{t***}:Y^{**}\times X^{**}\rightarrow Z^{**}$.
 The mapping $m^{t***t}:X^{**}\times Y^{**}\rightarrow Z^{**}$ in general is not equal to $m^{***}$, see [1], if $m^{***}=m^{t***t}$, then $m$ is called Arens regular. The mapping $y^{\prime\prime}\rightarrow m^{t***t}(x^{\prime\prime},y^{\prime\prime})$ is $weak^*-to-weak^*$ continuous for every $y^{\prime\prime}\in Y^{**}$, but the mapping $x^{\prime\prime}\rightarrow m^{t***t}(x^{\prime\prime},y^{\prime\prime})$ from $X^{**}$ into $Z^{**}$ is not in general  $weak^*-to-weak^*$ continuous for every $y^{\prime\prime}\in Y^{**}$. So we define the second topological center of $m$ as
$$Z_2(m)=\{y^{\prime\prime}\in Y^{**}:~~x^{\prime\prime}\rightarrow m^{t***t}(x^{\prime\prime},y^{\prime\prime})~~is~~weak^*-to-weak^*~~continuous\}.$$
It is clear that $m$ is Arens regular if and only if $Z_1(m)=X^{**}$ or $Z_2(m)=Y^{**}$. Arens regularity of $m$ is equivalent to the following
$$\lim_i\lim_j \langle  z^\prime,m(x_i,y_j)\rangle=\lim_j\lim_i \langle  z^\prime,m(x_i,y_j)\rangle,$$
whenever both limits exist for all bounded sequences $(x_i)_i\subseteq X$ , $(y_i)_i\subseteq Y$ and $z^\prime\in Z^*$, see [5].\\
The mapping $m$ is left strongly Arens irregular if $Z_1(m)=X$ and $m$ is right strongly Arens irregular if $Z_2(m)=Y$.\\
Let now $B$ be a Banach $A-bimodule$, and let\\
$$\pi_\ell:~A\times B\rightarrow B~~~and~~~\pi_r:~B\times A\rightarrow B.$$
be the left and right module actions of $A$ on $B$, respectively. Then $B^{**}$ is a Banach $A^{**}-bimodule$ with module actions
$$\pi_\ell^{***}:~A^{**}\times B^{**}\rightarrow B^{**}~~~and~~~\pi_r^{***}:~B^{**}\times A^{**}\rightarrow B^{**}.$$
Similarly, $B^{**}$ is a Banach $A^{**}-bimodule$ with module actions\\
$$\pi_\ell^{t***t}:~A^{**}\times B^{**}\rightarrow B^{**}~~~and~~~\pi_r^{t***t}:~B^{**}\times A^{**}\rightarrow B^{**}.$$
We may therefore define the topological centers of the left and right module actions of $A$ on $B$ as follows:\\
$$Z_{B^{**}}(A^{**})=Z(\pi_\ell)=\{a^{\prime\prime}\in A^{**}:~the~map~~b^{\prime\prime}\rightarrow \pi_\ell^{***}(a^{\prime\prime}, b^{\prime\prime})~:~B^{**}\rightarrow B^{**}$$$$~is~~~weak^*-to-weak^*~continuous\}$$
$$Z_{B^{**}}^t(A^{**})=Z(\pi_r^t)=\{a^{\prime\prime}\in A^{**}:~the~map~~b^{\prime\prime}\rightarrow \pi_r^{t***}(a^{\prime\prime}, b^{\prime\prime})~:~B^{**}\rightarrow B^{**}$$$$~is~~~weak^*-to-weak^*~continuous\}$$
$$Z_{A^{**}}(B^{**})=Z(\pi_r)=\{b^{\prime\prime}\in B^{**}:~the~map~~a^{\prime\prime}\rightarrow \pi_r^{***}(b^{\prime\prime}, a^{\prime\prime})~:~A^{**}\rightarrow B^{**}$$$$~is~~~weak^*-to-weak^*~continuous\}$$
$$Z_{A^{**}}^t(B^{**})=Z(\pi_\ell^t)=\{b^{\prime\prime}\in B^{**}:~the~map~~a^{\prime\prime}\rightarrow \pi_\ell^{t***}(b^{\prime\prime}, a^{\prime\prime})~:~A^{**}\rightarrow B^{**}$$$$~is~~~weak^*-to-weak^*~continuous\}$$

\noindent We note also that if $B$ is a left(resp. right) Banach $A-module$ and $\pi_\ell:~A\times B\rightarrow B$~(resp. $\pi_r:~B\times A\rightarrow B$) is left (resp. right) module action of $A$ on $B$, then $B^*$ is a right (resp. left) Banach $A-module$. \\
We write $ab=\pi_\ell(a,b)$, $ba=\pi_r(b,a)$, $\pi_\ell(a_1a_2,b)=\pi_\ell(a_1,a_2b)$, \\ $\pi_r(b,a_1a_2)=\pi_r(ba_1,a_2)$,~
$\pi_\ell^*(a_1b^\prime, a_2)=\pi_\ell^*(b^\prime, a_2a_1)$,~
$\pi_r^*(b^\prime a, b)=\pi_r^*(b^\prime, ab)$,~ \\for all $a_1,a_2, a\in A$, $b\in B$ and  $b^\prime\in B^*$
when there is no confusion.\\
Regarding $A$ as a Banach $A-bimodule$, the operation $\pi:A\times A\rightarrow A$ extends to $\pi^{***}$ and $\pi^{t***t}$ defined on $A^{**}\times A^{**}$. These extensions are known, respectively, as the first(left) and the second (right) Arens products, and with each of them, the second dual space $A^{**}$ becomes a Banach algebra. In this situation, we shall also simplify our notations. So the first (left) Arens product of $a^{\prime\prime},b^{\prime\prime}\in A^{**}$ shall be simply indicated by $a^{\prime\prime}b^{\prime\prime}$ and defined by the three steps:
 $$ \langle  a^\prime a,b\rangle= \langle  a^\prime ,ab\rangle,$$
  $$ \langle  a^{\prime\prime} a^\prime,a\rangle= \langle  a^{\prime\prime}, a^\prime a\rangle,$$
  $$ \langle  a^{\prime\prime}b^{\prime\prime},a^\prime\rangle= \langle  a^{\prime\prime},b^{\prime\prime}a^\prime\rangle.$$
 for every $a,b\in A$ and $a^\prime\in A^*$. Similarly, the second (right) Arens product of $a^{\prime\prime},b^{\prime\prime}\in A^{**}$ shall be  indicated by $a^{\prime\prime}ob^{\prime\prime}$ and defined by :
 $$ \langle  a oa^\prime ,b\rangle= \langle  a^\prime ,ba\rangle,$$
  $$ \langle  a^\prime oa^{\prime\prime} ,a\rangle= \langle  a^{\prime\prime},a oa^\prime \rangle,$$
  $$ \langle  a^{\prime\prime}ob^{\prime\prime},a^\prime\rangle= \langle  b^{\prime\prime},a^\prime ob^{\prime\prime}\rangle.$$
  for all $a,b\in A$ and $a^\prime\in A^*$.\\
  The regularity of a normed algebra $A$ is defined to be the regularity of its algebra multiplication when considered as a bilinear mapping. Let $a^{\prime\prime}$ and $b^{\prime\prime}$ be elements of $A^{**}$, the second dual of $A$. By $Goldstine^,s$ Theorem [4, P.424-425], there are nets $(a_{\alpha})_{\alpha}$ and $(b_{\beta})_{\beta}$ in $A$ such that $a^{\prime\prime}=weak^*-\lim_{\alpha}a_{\alpha}$ ~and~  $b^{\prime\prime}=weak^*-\lim_{\beta}b_{\beta}$. So it is easy to see that for all $a^\prime\in A^*$,
$$\lim_{\alpha}\lim_{\beta} \langle  a^\prime,\pi(a_{\alpha},b_{\beta})\rangle= \langle  a^{\prime\prime}b^{\prime\prime},a^\prime\rangle$$ and
$$\lim_{\beta}\lim_{\alpha} \langle  a^\prime,\pi(a_{\alpha},b_{\beta})\rangle= \langle  a^{\prime\prime}ob^{\prime\prime},a^\prime\rangle,$$
where $a^{\prime\prime}b^{\prime\prime}$ and $a^{\prime\prime}ob^{\prime\prime}$ are the first and second Arens products of $A^{**}$, respectively, see [16, 20].\\
  We find the usual first and second topological center of $A^{**}$, which are
  $$Z_1({A^{**}})=Z_{A^{**}}(A^{**})=Z(\pi)=\{a^{\prime\prime}\in A^{**}: b^{\prime\prime}\rightarrow a^{\prime\prime}b^{\prime\prime}~ is~weak^*-to-weak^*$$$$~continuous\},$$
   $$Z_2({A^{**}})=Z^t_{A^{**}}(A^{**})=Z(\pi^t)=\{a^{\prime\prime}\in A^{**}: a^{\prime\prime}\rightarrow a^{\prime\prime}ob^{\prime\prime}~ is~weak^*-to-weak^*$$$$~continuous\}.$$
Recall that  a left approximate identity $(=LAI)$ [resp. right
approximate identity $(=RAI)$]
in Banach algebra $A$ is a net $(e_{\alpha})_{{\alpha}\in I}$ in $A$ such that   $e_{\alpha}a\longrightarrow a$ [resp. $ae_{\alpha}\longrightarrow a$]. We
say that a  net $(e_{\alpha})_{{\alpha}\in I}\subseteq A$ is a
approximate identity $(=AI)$ for $A$ if it is $LAI$ and $RAI$ for $A$. If $(e_{\alpha})_{{\alpha}\in I}$ in $A$ is bounded and $AI$ for $A$, then we say that $(e_{\alpha})_{{\alpha}\in I}$ is a bounded approximate identity  ($=BAI$) for $A$.  Let $A$ have a $BAI$. If the
equality $A^*A=A^*,~~(AA^*=A^*)$ holds, then we say that $A^*$
factors on the left (right). If both equalities $A^*A=AA^*=A^*$
hold, then we say
that $A^*$  factors on both sides.\\An element $e^{\prime\prime}$ of $A^{**}$ is said to be a mixed unit if $e^{\prime\prime}$ is a
right unit for the first Arens multiplication and a left unit for
the second Arens multiplication. That is, $e^{\prime\prime}$ is a mixed unit if
and only if,
for each $a^{\prime\prime}\in A^{**}$, $a^{\prime\prime}e^{\prime\prime}=e^{\prime\prime}o a^{\prime\prime}=a^{\prime\prime}$. By
[4, p.146], an element $e^{\prime\prime}$ of $A^{**}$  is  mixed
      unit if and only if it is a $weak^*$ cluster point of some BAI $(e_\alpha)_{\alpha \in I}$  in
      $A$.\\

\begin{center}
\section{ \bf Factorization property and  topological centers of module actions }
\end{center}
Baker, Lau and Pym in [3] proved that for Banach algebra $A$ with bounded right approximate identity, $(A^*A)^\bot$ is an ideal of right annihilators in $A^{**}$ and $A^{**}\cong (A^*A)\oplus(A^*A)^\bot$.
In the following, for a Banach $A-bimodule$ $B$, we  study  the similar discussion on the module actions and for Banach $A-bimodule$ $B$,  we show that $$B^{**}=(B^*A)^*\oplus(B^*A)^\bot.$$\\

\noindent{\it{\bf Theorem 2-1.}} Let $B$ be a Banach $A-bimodule$ and $A$ has a $BRAI$. Then the following assertions are hold:\\
i) $(B^*A)^\bot=\{b^{\prime\prime}\in B^{**}:~\pi_\ell^{***}(a^{\prime\prime},b^{\prime\prime})=0~~for~all~a^{\prime\prime}\in A^{**}\}$.\\
ii) $(B^*A)^*$ is isomorphism with $Hom_A(B^*,A^*)$.
\begin{proof} i) Let $b^{\prime\prime}\in(B^*A)^\bot$. Then for all $b^\prime\in B^*$ and $a\in A$, we have
$$ \langle  \pi^{**}_\ell(b^{\prime\prime},b^\prime),a\rangle= \langle  b^{\prime\prime},\pi^{*}_\ell(b^\prime,a)\rangle=0,$$
it follows that for all $a^{\prime\prime}\in A^{**}$, we have
$$ \langle  \pi^{***}_\ell(a^{\prime\prime},b^{\prime\prime}),b^\prime\rangle=
 \langle  a^{\prime\prime},\pi^{**}_\ell(b^{\prime\prime},b^\prime)\rangle=0.$$
Conversely, let $b^{\prime\prime}\in B^{**}$ such that $\pi^{***}_\ell(a^{\prime\prime},b^{\prime\prime})=0$ for all
$a^{\prime\prime}\in A^{**}$. Then for all $a\in A$ and $b^\prime\in B^*$, we have
$$ \langle  b^{\prime\prime},\pi^*_\ell(b^\prime,a)\rangle= \langle  \pi^{**}_\ell(b^{\prime\prime},b^\prime),a\rangle
= \langle  a,\pi^{**}_\ell(b^{\prime\prime},b^\prime)\rangle
=
 \langle  \pi^{***}_\ell(a^{},b^{\prime\prime}),b^\prime\rangle=0,$$
which implies that $b^{\prime\prime}\in(B^*A)^\bot$.\\
ii) Suppose that $b^{\prime\prime}\in B^{**}$. We define $T_{b^{\prime\prime}}\in Hom_A(B^*,A^*)$, that is,
$T_{b^{\prime\prime}}b^\prime=\pi^{**}_\ell(b^{\prime\prime},b^\prime)$. Then $\Lambda=b^{\prime\prime}\rightarrow T_{b^{\prime\prime}}$ is linear continuous map from $B^{**}$ into $Hom_A(B^*,A^*)$ such that
$$Ker \Lambda=\{b^{\prime\prime}\in B^{**}:~\pi^{**}_\ell(b^\prime,b^{\prime\prime})=0~~for~all~b^\prime\in B^*\}.$$
Consequently, $b^{\prime\prime}\in Ker \Lambda$ if and only if
$$ \langle  b^{\prime\prime},\pi^*_\ell(b^\prime,a)= \langle  \pi^{**}_\ell(b^{\prime\prime},b^\prime),a\rangle=0,$$
where $b^\prime\in B^*$ and $a\in A$. It follows that $b^{\prime\prime}\in(B^*A)^\bot$. Since $(B^*A)^*\cong\frac{B^{**}}{(B^*A)^\bot}$, the continuous linear mapping $\Lambda$ from $(B^*A)^*$ into $Hom_A(B^*,A^*)$ is injective.\\
Conversely, suppose that $T\in Hom_A(B^*,A^*)$ and $e^{\prime\prime}\in A^{**}$ is any right identity for $A^{**}$. We define $b^{\prime\prime}_T\in B^{**}$ such that for all $b^\prime$, we have
$$ \langle  b^{\prime\prime}_T,b^\prime\rangle= \langle  e^{\prime\prime},Tb^\prime\rangle.$$
It is clear that the linear mapping $T\rightarrow b^{\prime\prime}_T$ is continuous. For all $a\in A$, we have
$$ \langle  \pi^{**}_\ell(b^{\prime\prime}_T,b^\prime),a\rangle= \langle  b^{\prime\prime}_T,\pi^{*}_\ell(b^\prime,a)\rangle=
 \langle  e^{\prime\prime},T\pi^*_\ell(b^\prime,a)\rangle= \langle  e^{\prime\prime},(Tb^\prime)a\rangle$$
$$= \langle  ae^{\prime\prime},Tb^\prime\rangle= \langle  Tb^\prime,a\rangle.$$
Consequently, $\pi^{**}_\ell(b^{\prime\prime}_T,b^\prime)=Tb^\prime$. It follows that the linear mapping  $T\rightarrow b^{\prime\prime}_T\rightarrow T_{b^{\prime\prime}_T}$ is the identity map and consequently the isomorphism between $Hom_A(B^*,A^*)$ and $(B^*A)^*$ is established.\\
\end{proof}
\noindent By using  proceeding theorem we observe that $Hom_A(B^*,A^*)$  has an right identity. Also if $B=A$, we obtain Theorem 1-1 from [3].\\\\

\noindent{\it{\bf Corollary 2-2.}} Let $B$ be a Banach $A-bimodule$ and let $e^{\prime\prime}$ be any right identity of $A^{**}$. Then $e^{\prime\prime}B^{**}\cong (B^*A)^*$ and  $(B^*A)^\bot=\{b^{\prime\prime}-e^{\prime\prime}b^{\prime\prime}:~
b^{\prime\prime}\in B^{**}\}$. Thus $B^{**}= (B^*A)^*\oplus  (B^*A)^\bot$.\\\\

\noindent{\it{\bf Example 2-3.}}  Let $G$ be a locally compact group. Let $1\leq p < \infty$ and $\frac{1}{p}+\frac{1}{q}=1$. Then by using Theorem 2-1, we conclude that
$$(L^p(G)*L^1(G))^\bot=\{b\in L^q(G):~a^{\prime\prime}b=0~for ~every~a^{\prime\prime}\in L^\infty (G)\},$$
$$(L^p(G)*L^1(G))^*\cong Hom_{L^1(G)} (L^p(G),L^\infty(G)).$$\\\\

\noindent{\it{\bf Theorem 2-4.}} Assume that $B$ is a left Banach  $A-module$ and $A$ has a $BAI$. Then  we have the following assertions.
\begin{enumerate}

\item  If $B^*$ factors on the left and $B^{**}$ has a left unit as $A^{**}-module$, then $(B^*)^\perp=0$.
\item  If $B^*$ not factors on the left  and $e^{\prime\prime}$ is a left unit as $A^{**}-module$ in $B^{**}$, then  $e^{\prime\prime}\notin Z^\ell_{B^{**}}(A^{**})$.\\

\end{enumerate}

\begin{proof}
\begin{enumerate}

\item  Let $a\in A$,  $b^\prime\in B^*$ and $b^{\prime\prime}\in (B^*A)^\perp$. Then
$$ \langle  \pi_\ell^{**}(b^{\prime\prime},b^\prime),a\rangle= \langle  b^{\prime\prime},\pi_\ell^{*}(b^\prime,a)\rangle=0.$$
Thus, for all $b^{\prime\prime}\in A^{**}$, we have
$$ \langle  \pi_\ell^{***}(a^{\prime\prime},b^{\prime\prime}),b^\prime\rangle=
 \langle  a^{\prime\prime},\pi_\ell^{**}(b^{\prime\prime},b^\prime)\rangle=0.$$
It follows that $\pi_\ell^{***}(a^{\prime\prime},b^{\prime\prime})=0$.\\
 Now let $e^{\prime\prime}\in A^{**}$ be a left unit as $A^{**}-module$ for $B^{**}$, then  we have
$$b^{\prime\prime}=\pi_\ell^{***}(e^{\prime\prime},b^{\prime\prime})=0.$$
Thus proof  hold.
\item  Assume a contradiction that $e^{\prime\prime}\in Z_{B^{**}}(A^{**})$. Take $b^{\prime\prime}\in ({B^*A})^\perp$. Then
    $$b^{\prime\prime}=\pi_\ell^{***}(e^{\prime\prime},b^{\prime\prime})
    =\pi_\ell^{t***t}(e^{\prime\prime},b^{\prime\prime}).$$
Let $(e_\alpha)_\alpha\subseteq A$ such that $e_\alpha\stackrel{w^*} {\rightarrow}e^{\prime\prime}$ and let $b\in B$, $b^\prime\in B^*$. Then
$$  \langle  \pi_\ell^{t**}(e^{\prime\prime},b^\prime),b\rangle= \langle  e^{\prime\prime},\pi_\ell^{t*}(b^\prime,b)\rangle=
lim_\alpha \langle  \pi_\ell^{t*}(b^\prime,b),e_\alpha\rangle$$
$$=lim_\alpha \langle  b^\prime,\pi_\ell^{t}(b,e_\alpha)\rangle=lim_\alpha \langle  b^\prime,\pi_\ell^{}(e_\alpha,b)\rangle
=lim_\alpha \langle  \pi_\ell^{*}(b^\prime,e_\alpha),b\rangle$$
It follows that
$$w^*-lim_\alpha\pi_\ell^{*}(b^\prime,e_\alpha)=\pi_\ell^{t**}(e^{\prime\prime},b^\prime).$$
Let $(b_\beta)_\beta\subseteq B$ such that $b_\beta\stackrel{w^*} {\rightarrow}b^{\prime\prime}$. Then
$$  \langle  b^{\prime\prime},b^\prime\rangle= \langle  \pi_\ell^{t***t}(e^{\prime\prime},b^{\prime\prime}),b^\prime\rangle= \langle  \pi_\ell^{t***}(b^{\prime\prime},e^{\prime\prime}),b^\prime\rangle=
 \langle  b^{\prime\prime},\pi_\ell^{t**}(e^{\prime\prime},b^\prime)\rangle$$
$$=lim_\beta \langle  \pi_\ell^{t**}(e^{\prime\prime},b^\prime),b_\beta\rangle=
lim_\beta lim_\alpha \langle  \pi_\ell^{*}(b^\prime,e_\alpha),b_\beta\rangle$$
$$=lim_\beta lim_\alpha \langle  b^\prime,\pi_\ell^{}(e_\alpha,b_\beta)\rangle= lim_\alpha lim_\beta \langle  b^\prime,\pi_\ell^{}(e_\alpha,b_\beta)\rangle$$
$$=lim_\alpha \langle  b^\prime,\pi_\ell^{***}(e_\alpha,b^{\prime\prime})\rangle
=lim_\alpha \langle  \pi_\ell^{***}(e_\alpha,b^{\prime\prime}),b^\prime\rangle$$
$$=lim_\alpha \langle  e_\alpha,\pi_\ell^{**}(b^{\prime\prime},b^\prime)\rangle
=lim_\alpha \langle  b^{\prime\prime},\pi_\ell^{*}(b^\prime,e_\alpha)\rangle=0.$$
It follows that $({B^*A})^\perp=0$. By using Corollary 2-2, we have  $B^{**}=(B^*A)^\perp\oplus ({B^*A})^*$, and so  $B^{**}=({B^*A})^*$. Now since ${B^*A}\neq B^*$, by Hahn Banach theorem, there is $0\neq b^{\prime\prime}\in B^{**}$ such that $b^{\prime\prime}|_{B^*A}=0$. It follows that $b^{\prime\prime}\in (B^*A)^\perp$ which  is contradiction.\\
\end{enumerate}
\end{proof}
\noindent{\it{\bf Corollary 2-5.}} For a left Banach $A-module$ $B$, we have the following statements.
\begin{enumerate}

\item  If $\overline{B^*A}= B^*$ and $B^{**}$ has a left unit as $A^{**}-module$, then $(B^*)^\perp=0$.
\item  If $\overline{B^*A}\neq B^*$ and $e^{\prime\prime}$ is a left unit as $A^{**}-module$ in $B^{**}$, then  $e^{\prime\prime}\notin Z^\ell_{B^{**}}(A^{**})$.
\end{enumerate}
Proof is similar to  Theorem 2-4.\\\\

\noindent{\it{\bf Corollary 2-6.}} Assume that $B$ is a left Banach  $A-module$ and  $ B^{**}$ has a left unit $A^{**}-module$. If $\overline{B^*A}\neq B^*$, then $Z^\ell_{B^{**}}(A^{**})\neq A^{**}$.\\\\

In the proceeding corollary,  take $B=A$ and assume that $A^*$ not factors on the left. Then, if $A$ has a $LBAI$, we conclude that $A$ is not Arens regular.\\\\

\noindent{\it{\bf Example 2-7.}}   Let $G$ be a  locally compact group, by using Corollary 2-6, we have the following inequality
~~~~~~~~~$$Z^\ell_{L^1(G)^{**}}(M(G)^{**})\neq M(G)^{**} ~~,~~~~~ Z^\ell_{M(G)^{**}}(L^1(G)^{**})\neq L^1(G)^{**}.$$
Thus  $L^1(G)$ and $M(G)$ are not Arens regular.\\\\

Baker,  Dales, Lau, Pym and \"{U}lger in [3, 5, 6, 16] have studied topological center of some Banach algebras and they have solved some problems related to module homomorphism and topological centers of Banach algebras. In the following we extend some problems of them into module actions with some new result  and applications in  group algebras.\\
Throughout this paper, the notations $WSC$ is used for weakly sequentially
complete Banach space $A$, that is, $A$ is said to be weakly
sequentially complete, if every weakly Cauchy sequence in $A$ has
a weak limit in $A$.\\
Assume that  $B$ is a Banach $A-bimodule$. We say that  $B$ factors on the left (right) with respect to $A$ if $B=BA~(B=AB)$. We say that $B$ factors on both sides, if $B=BA=AB$.\\
\noindent Suppose that $A$ is a Banach algebra and $B$ is a Banach $A-bimodule$. According to [28] $B^{**}$ is a Banach $A^{**}-bimodule$, where  $A^{**}$ is equipped with the first Arens product.
 We define  $B^{*}B$ as a subspace of $A$, that is, for all $b^{\prime}\in B^{*}$ and $b\in B$, we define
$$ \langle  b^{\prime}b,a\rangle= \langle  b^{{\prime}},ba\rangle;$$
We similarly define $B^{***}B^{**}$ as a subspace of $A^{**}$  and we take $A^{(0)}=A$ and $B^{(0)}=B$.\\
Let $B$ be a left Banach $A-module$ and $(e_{\alpha})_{\alpha}\subseteq A$ be a LAI [resp. weakly left approximate identity(=WLAI)] for $A$. We say that $(e_{\alpha})_{\alpha}$ is left approximate identity  ($=LAI$)[ resp. weakly left approximate identity  (=$WLAI$)] for $B$, if for all $b\in B$, we have $\pi_\ell (e_{\alpha},b) \stackrel{} {\rightarrow}
b$ ( resp. $\pi_\ell (e_{\alpha},b) \stackrel{w} {\rightarrow}
b$). The definition of the right approximate identity ($=RAI$)[ resp. weakly right approximate identity ($=WRAI$)] is similar.\\
We say that $(e_{\alpha})_{\alpha}$ is a approximate identity  ($=AI$)[ resp. weakly approximate identity  ($WAI$)] for $B$, if $B$ has left and right approximate identity  [ resp. weakly left and right approximate identity ] that are equal.\\
Let $B$ be a left Banach as $A-module$ and  $e$ be a left  unit element of $A$. Then we say that $e$ is a left unit (resp. weakly left unit)  as $A-module$ for $B$, if $\pi_\ell(e,b)=b$ (resp. $ \langle  b^\prime , \pi_\ell(e,b)\rangle= \langle  b^\prime , b\rangle$ for all $b^\prime\in B^*$) where $b\in B$. The definition of right unit (resp. weakly right unit) as $A-module$ is similar.\\
We say that a Banach $A-bimodule$ $B$ is an unital, if $B$ has the same left and right unit as $A-module$.\\\\

\noindent{\it{\bf Lemma 2-8.}} Let $B$ be a   Banach $A-bimodule$. Suppose that $A$ has a $BAI$ $(e_\alpha)_\alpha\subseteq A$. Then
\begin{enumerate}
\item $B$ factors on the left if and only if $\pi_r(b,e_\alpha) \stackrel{w} {\rightarrow}b$ for every $b\in B$.
\item $B$ factors on the right   if and only if $\pi_\ell(e_\alpha, b) \stackrel{w} {\rightarrow}b$ for every $b\in B$.
\item  If $B^*$ factors on the right, then  $\pi_r(b,e_\alpha) \stackrel{w} {\rightarrow}b$ for every $b\in B$.

\end{enumerate}
\begin{proof}
\begin{enumerate}
\item  Suppose that $B$ factors on the left. Then for every   $b\in B$, there are $y\in B$ and $a\in A$ such that $b=ya$. Thus for every $b^\prime \in B^*$, we have
$$ \langle  b^\prime ,\pi_r(b,e_\alpha)\rangle= \langle  b^\prime ,\pi_r(ya,e_\alpha)\rangle= \langle  b^\prime ,\pi_r(y,ae_\alpha)\rangle= \langle  \pi^*_r(b^\prime ,y),ae_\alpha\rangle$$$$\rightarrow
 \langle  \pi^*_r(b^\prime ,y),a\rangle= \langle  b^\prime ,ya\rangle= \langle  b^\prime ,b\rangle.$$
It follows that $\pi_r(b,e_\alpha)\stackrel{w} {\rightarrow}b$. \\
Conversely, by $Cohen^,s$ factorization Theorem, since  $BA$ is a closed subspace of $B$, the proof is hold.
\item Proof is similar to (1).
\item Assume that $B^*$ factors on the right with respect to $A$. Then for every   $b^\prime \in B^*$, there are $y^\prime \in B$ and $a\in A$ such that $b^\prime =ay^\prime $.
Consequently for every $b\in B$, we have
$$ \langle  b^\prime ,\pi_r(b,e_\alpha)\rangle= \langle  ay^\prime ,\pi_r(b,e_\alpha)\rangle= \langle  y^\prime ,\pi_r(b,e_\alpha) a\rangle= \langle  y^\prime ,\pi_r(b,e_\alpha a) \rangle  $$$$=
 \langle  \pi^*_r(y^\prime ,b),e_\alpha a\rangle\rightarrow   \langle  \pi^*_r(y^\prime ,b), a\rangle=\langle  y^\prime, \pi_r(b, a)\rangle= \langle  ay^\prime ,b \rangle$$$$=
 \langle  b^\prime ,b \rangle.$$
It follows that $\pi_r(b,e_\alpha) \stackrel{w} {\rightarrow}b$. \\

\end{enumerate}
\end{proof}

In the proceeding theorem, if we take $B=A$, then we obtain Lemma 2.1 from [16].\\\\

\noindent {\it{\bf Theorem 2-9.}} Let  $B$ be a  Banach  $A-bimodule$ and $A$ has a sequential $WBAI$. Then we have the following assertions.\\
(i) Let $B^*$ be a $WSC$ and $A^*$ factors on the left. Then
\begin{enumerate}
\item  If  $B$ factors on the right, it follows that $B^*$ factors on the left.
\item If  $B^*$ factors on the right, it follows that $B$ factors on the left.
\end{enumerate}
(ii) Let $B^{**}B^*=A^{**}A^*$. Then $A^*$ factors on the left if and only if $B^*$ factors on the left.\\
(iii) Suppose that $A$ is $WSC$ and $B$ factors on the left (resp. right). If $B^*B=A^*$, then we have the following assertions.

\begin{enumerate}
\item   $A$ is an unital and  $B$ has a right (resp. left) unit as  Banach $A-module$.
\item     $A^{*}$ factors on the both side and  $B^*$  factors  on the  right (resp. left).
\item   $B^{**}\cong(AB^*)^*$ (resp. $B^{**}\cong(B^*A)^*$).

\end{enumerate}
\begin{proof}
i) 1)  Assume that $b^{\prime\prime}\in B^{**}$ and $b^\prime\in B^*$. Since $A^*$ factors on the left, there are $a^\prime\in A^*$ and $a\in A$ such that $b^{\prime\prime}b^\prime=a^{\prime}a$. Suppose that $(e_n)_n\subseteq A$ is a sequential $WBAI$ for $A$. Then we have
$$ \langle  b^{\prime\prime},b^\prime e_n\rangle= \langle  b^{\prime\prime}b^\prime ,e_n\rangle= \langle  a^{\prime}a ,e_n\rangle= \langle  a^{\prime},a e_n\rangle\rightarrow
 \langle  a^{\prime},a \rangle.$$
It follows that  the sequence $(b^\prime e_n)_n$ is weakly Cauchy sequence in $B^*$. Since $B^*$ is $WSC$, there is $x^\prime\in B^*$ such that $b^\prime e_n \stackrel{w} {\rightarrow}x^\prime$. On the other hand,  since $B$ factors on the right, by using Lemma 2-8, for each $b\in B$, we have $e_n b \stackrel{w} {\rightarrow}b$. Then we have
$$ \langle  x^\prime ,b\rangle= \lim_n\langle  b^\prime e_n ,b\rangle=  \lim_n\langle  b^\prime ,e_n b\rangle=  \langle  b^\prime ,b\rangle.$$
It follows that $x^\prime=b^\prime$, and so by Lemma 2-8, $B^*$ factors on the left.\\
i) 2) Proof is similar to part (1).\\
ii) Let $a^{\prime\prime}\in A^{**}$ and $a^\prime\in A^*$. Then there are $b^{\prime\prime}\in B^{**}$ and $b^\prime\in B^*$ such that $b^{\prime\prime}b^\prime=a^{\prime\prime}a^\prime$. Then
$$ \langle  a^{\prime\prime},a^\prime e_n\rangle= \langle  a^{\prime\prime}a^\prime, e_n\rangle= \langle  b^{\prime\prime}b^\prime, e_n\rangle= \langle  b^{\prime\prime},b^\prime e_n\rangle.$$
Thus, by $Cohen^,s$ factorization theorem, proof hold.\\
iii)  1) Suppose that $(e_k)_k\subseteq A$ is a sequential $WBAI$ for $A$. Let $a^\prime \in A^*$. Since $B^*B=A^*$, there are  $b^\prime\in B^*$ and $b\in B$ such that $b^\prime b=a^\prime$. Since  $B$ factors on the left,   there are $y\in B$ and $a\in A$ such that $b=ya$. Then we have
$$ \langle  a^\prime,e_k\rangle= \langle  b^\prime b,e_k\rangle= \langle  b^\prime ,be_k\rangle= \langle  b^\prime ,yae_k\rangle$$$$=
 \langle  b^\prime y,ae_k\rangle\rightarrow \langle  b^\prime y,a\rangle= \langle  b^\prime ,ya\rangle
$$$$= \langle  b^\prime ,b\rangle.$$
This shows that the sequence  $(e_k)_k\subseteq A$ is weakly sequence in $A$. Since $A$ is WSC, it convergence weakly to some element $e$ of $A$. Then, for each $x\in A$, we have
$$xe=x(w-\lim_k e_k)=w-\lim_kxe_k=a.$$
It is similar to that $ex=x$, and so $A$ is unital.\\
Now let $b\in B$, then

$$ \langle  b^\prime , be\rangle= \langle  b^\prime  b,e\rangle=\lim_k \langle  b^\prime b, e_k\rangle=\lim_k \langle  b^\prime ,ya e_k\rangle$$$$= \lim_k \langle  b^\prime y,a e_k\rangle\rightarrow  \langle  b^\prime y,a \rangle= \langle  b^\prime ,b \rangle.$$
Thus $be=b$ for all $b\in B$.

iii)  2)  By using part (1) and [16, Theorem 2.6], it is clear that $A^{*}$ factors on the both side. Now let $b^\prime\in B^*$ and $b\in B$. By part (1), set $e\in A$ as a left unite element of $B$. Then
$$ \langle  eb^\prime , b\rangle= \langle  b^\prime , be\rangle= \langle  b^\prime , b\rangle.$$
It follows that  $eb^\prime =b^\prime$. Thus $B^*$ factors on the right.\\
  iii)  3) Now let $b^{\prime\prime}\in (AB^*)^\perp$. By using part (2), since $B^*$ factors on the right, for every $b^\prime \in B^*$ there are $x^\prime \in B^*$ and  $a\in A$ such that
$b^\prime=ax^\prime $. Then
$$ \langle  b^{\prime\prime}, b^\prime\rangle= \langle  b^{\prime\prime}, ax^\prime \rangle=0.$$
It follows that $b^{\prime\prime}=0$. It follows that $(AB^*)^\perp=\{0\}$. Therefore by using Corollary 2-2, we are done.\\\end{proof}

In the proceeding theorem if we take $A=B$, then we have the following statements:\\
Let   $A$ has a sequential $WBAI$ and Suppose that  $A$ is a $WSC$. Then
\begin{enumerate}
\item   $A^*$  factors on the left (resp. right) if and only if $A^*$ factors on the right (resp. left), see [16].
\item   $A^*$ factors on the left (resp. right)  if and only if   $A$ is an unital, see [16].
\item  If  $A^*$ factors on the left (resp. right), then $(AA^*)^\perp\cong A^{**}$ (resp. $(A^*A)^\perp\cong A^{**}$).\\\\
\end{enumerate}

\noindent {\it{\bf Example 2-10.}}\\
i) Let $G$ be a locally compact groups. Take $B=c_0(G)$ and $B^*=A=\ell^1(G)$. Since $\ell^1(G)\ell^1(G)=\ell^1(G)$, by proceeding theorem we have  $c_0(G)\ell^1(G)=c_0(G)$.\\
ii) Let $\omega:G\rightarrow \mathbb{R}^+\setminus\{0\}$ be a continuous function on a locally compact group $G$. Suppose that
$$L^1(G,\omega)=\{f ~~is~~Borel~~measurable:~~\parallel f\parallel_\omega=\int_G \mid f(s)\omega(s)\mid ds <  \infty\},~see~[6].$$
Let $X$ be a closed separable subalgebra of $L^1(G,\omega)$. Then by using [6, Theorem 7.1] and [24, Lemma 3.2], X is $WSC$ and it has a sequential $BAI$, respectively. Therefore by above results if $X^*X=X^*$, then  $X$ is unital and $(XX^*)^\perp=\{0\}$. On the other hand, if the identity of $G$ has a countable neighborhood base, then $L^1(G,\omega)$ has a sequential $BAI$. Now let $L^1(G,\omega)^*=L^\infty(G,\frac {1}{\omega})$ factors on the left. Hence
$$LUC(G,\frac {1}{\omega})=L^\infty(G,\frac {1}{\omega})L^1(G,\omega)=L^\infty(G,\frac {1}{\omega}),$$
where $L^\infty(G,\frac {1}{\omega})=\{f ~~is~~Borel~~measurable:~~\parallel f\parallel_{\infty ,\omega}=ess~sup_{s\in G} \frac{\mid f(s)\mid}{\omega(s)} <  \infty\}.$ Consequently by above results, $L^1(G,\omega)$ is unital and  $RUC(G,\frac {1}{\omega})^\bot=\{0\}$, where $$RUC(G,\frac {1}{\omega})=L^1(G,\omega) L^\infty(G,\frac {1}{\omega}).$$ Thus $L^1(G,\omega)^{**}=RUC(G,\frac {1}{\omega})^*$.\\\\

\noindent {\it{\bf Theorem 2-11.}} Suppose that $B$ is a left Banach $A-module$ and it has a $WLBAI$ $(e_\alpha)_\alpha\subseteq A$. Then  we have the following assertions.
\begin{enumerate}

\item $B$ factors on the left.
\item  If $A^*$ factors on the left, then $B^*$ factors on the left.

\end{enumerate}

\begin{proof}
\begin{enumerate}
\item  By using Lemma 2-8, proof hold.
\item Let $b^{\prime\prime}\in B^{**}$ and $b^\prime\in B^*$. Since $\pi_\ell^{**}(b^{\prime\prime},b^\prime )\in A^*$ and $A^*$ factors on the right, there are $a^\prime\in A^*$ and $a\in A$ such that $\pi_\ell^{**}(b^{\prime\prime},b^\prime )=a^\prime a$. Without loss generality,  we let $e_\alpha\stackrel{w^*} {\rightarrow}e^{\prime \prime}$ where $e^{\prime \prime}$ left unit for $A^{**}$. Then for every $b\in B$, we have
$$ \langle  \pi_\ell^{****}(b^\prime, e^{\prime \prime}),b\rangle=  \langle  b^\prime, \pi_\ell^{***}(e^{\prime \prime},b)\rangle=\lim_\alpha \langle  b^\prime,\pi_\ell^{}(e_\alpha b)\rangle= \langle  b^\prime,b\rangle.$$
It follows that $\pi_\ell^{****}(b^\prime, e^{\prime \prime})=b^\prime$. Now, we have the following equality
$$ \langle  b^{\prime \prime},\pi_\ell^{*}(b^\prime ,e_\alpha)-b^\prime\rangle= \langle  b^{\prime \prime},\pi_\ell^{****}(b^\prime ,(e_\alpha-e^{\prime \prime}))\rangle$$$$=
 \langle \pi_\ell^{*****}( b^{\prime \prime},b^\prime) ,(e_\alpha-e^{\prime \prime})\rangle= \langle \pi_\ell^{**}( b^{\prime \prime},b^\prime) ,(e_\alpha-e^{\prime \prime})\rangle$$$$= \langle  a^\prime a,(e_\alpha-e^{\prime \prime})\rangle=
 \langle  a^\prime ,ae_\alpha-ae^{\prime \prime}\rangle$$$$= \langle  a^\prime ,ae_\alpha-a\rangle\rightarrow 0.$$
It follows that $\pi_\ell^{*}(b^\prime, e_\alpha)\stackrel{w} {\rightarrow}b^\prime$, and so by $Cohen^,s$ factorization, we are done.\\ \end{enumerate} \end{proof}

\noindent {\it{\bf Theorem 2-12.}} Suppose that $B$ is a right Banach $A-module$ and it has a  $RBAI$ $(e_\alpha)_\alpha\subseteq A$. Then  we have the following assertions.
\begin{enumerate}

\item $B$ factors on the right.
\item  Let ${Z}^\ell_{A^{**}}(B^{**})=B^{**}$. If $A^*$ factors on the right, then $B^*$ factors on the right.
\end{enumerate}

\begin{proof}
\begin{enumerate}
\item By using Lemma 2-8, proof hold.
\item Let $b^{\prime\prime}\in B^{**}$ and $b^\prime\in B^*$. First we show that $\pi _r ^{****}(b^\prime, b^{\prime\prime})\in A^*$. Suppose that $(a^{\prime\prime}_\alpha)_\alpha \subseteq A^{**}$ such that $a^{\prime\prime}_\alpha \stackrel{w^*} {\rightarrow}a^{\prime \prime}$. Since ${Z}^\ell_{A^{**}}(B^{**})=B^{**}$, for each $b^{\prime \prime}\in B^{**}$, we have $\pi _r ^{***}(b^{\prime\prime},a^{\prime\prime}_\alpha )\stackrel{w^*} {\rightarrow}\pi _r ^{***}(b^{\prime\prime},a^{\prime \prime})$. Then
    $$ \langle  \pi _r ^{****}(b^\prime, b^{\prime\prime}), a_\alpha^{\prime \prime}\rangle=
 \langle   \pi _r ^{***}(b^{\prime\prime},a^{\prime\prime}_\alpha ),b^\prime\rangle\rightarrow  \langle   \pi _r ^{***}(b^{\prime\prime},a^{\prime\prime} ),b^\prime\rangle
= \langle  \pi _r ^{****}(b^\prime, b^{\prime\prime}), a^{\prime \prime}\rangle.$$
Consequently $\pi _r ^{****}(b^\prime, b^{\prime\prime})\in (A^{**},weak^*)^*=A^*$. Since $A^*$ factors on the right, there are $a^\prime\in A^*$ and $a\in A$ such that $\pi _r ^{****}(b^\prime, b^{\prime\prime})=a^\prime a$.
 Without loss generality,  we let $e_\alpha\stackrel{w^*} {\rightarrow}e^{\prime \prime}$ where $e^{\prime \prime}$ right unit for $A^{**}$. Then for each $b\in B$, we have
$$ \langle  \pi _r ^{**}(e^{\prime \prime},b^\prime) ,b\rangle= \langle  b^\prime ,\pi _r ^{***}(b,e^{\prime \prime})\rangle=\lim_\alpha \langle  b^\prime ,\pi_r(b,e_\alpha)\rangle= \langle  b^\prime, b\rangle.$$
It follows that $\pi _r ^{**}(e^{\prime \prime},b^\prime)=b^\prime$. Now we have the following equality
$$ \langle  b^{\prime \prime},\pi _r ^{**}(e_\alpha ,b^\prime)-b^\prime\rangle= \langle  b^{\prime \prime},\pi _r ^{**}(e_\alpha ,b^\prime)-\pi _r ^{**}(e^{\prime \prime},b^\prime)\rangle$$$$=
 \langle  \pi _r ^{****}(b^\prime, b^{\prime\prime}),e_\alpha-e^{\prime \prime}\rangle= \langle  a^\prime a,e_\alpha-e^{\prime \prime}\rangle$$$$=
 \langle  a^\prime ,ae_\alpha-ae^{\prime \prime}\rangle= \langle  a^\prime ,ae_\alpha-a\rangle\rightarrow 0.$$
Thus $\pi_r^{**}(e_\alpha, b^\prime)\stackrel{w} {\rightarrow}b^{\prime}$. Consequently by $Cohen^,s$ factorization, we are done.\\
\end{enumerate}
\end{proof}

\noindent {\it{\bf  Corollary 2-13.}} Suppose that $B$ is a  Banach $A-bimodule$ and it has a  $BAI$ $(e_\alpha)_\alpha\subseteq A$. Then  we have the following assertions.
\begin{enumerate}

\item $B$ factors.
\item  Let ${Z}^\ell_{A^{**}}(B^{**})=B^{**}$. If $A^*$ factors on the both side, then $B^*$ factors on the both side.\\\\

\end{enumerate}

\noindent {\it{\bf Example 2-14.}} Assume that $G$ is a locally compact group. We know that $L^1(G)$ is a $M(G)-bimodule$. Since $M(G)L^1(G)\neq M(G)$ and $L^1(G)M(G)\neq M(G)$, by using Theorem 2-11 and 2-12, we conclude that every $LBAI$ or $RBAI$ for $L^1(G)$ is not
$LBAI$ or $RBAI$ for $M(G)$, respectively. \\\\

\noindent {\it{\bf Definition 2-15.}} Let $B$ be a left Banach $A-module$. Then $B$ is said to be left weakly completely continuous $(=Lwcc)$, if for each $a\in A$, the mapping $b\rightarrow \pi_\ell(a,b)$ from $B$ into $B$ is weakly compact. The definition of right weakly completely continuous $(=Rwcc)$ is similar. We say that $B$ is a weakly completely continuous $(=wcc)$, if $B$ is $Lwcc$ and  $Rwcc$.\\\\

\noindent{\it{\bf Theorem 2-16.}} Let $B$ be a left Banach $A-module$ and for all $a\in A$,  $L_a$ be the linear mapping from $B$ into itself such that $L_ab=\pi_\ell(a,b)$ for all $b\in B$. Then $AB^{**}\subseteq B$ if and only if $L_a$ is weakly compact.

\begin{proof}
 Assume that $AB^{**}\subseteq B$. We take $L_a^*$ as  the adjoint of $L_a$. It is easy to show that $L^*_ab^\prime=\pi_\ell^*(b^\prime,a)$ for all $b^\prime\in B^*$. Then for every $b^{\prime\prime}\in B^{**}$, we have
$$
  \langle  L_a^{**}b^{\prime\prime},b^\prime\rangle= \langle  b^{\prime\prime},L_a^{*}b^\prime\rangle= \langle  b^{\prime\prime},\pi_\ell^{*}(b^\prime ,a)\rangle= \langle  \pi_\ell^{**}(b^{\prime\prime},b^\prime) ,a\rangle
$$
$$
= \langle  a,\pi_\ell^{**}(b^{\prime\prime},b^\prime)\rangle= \langle  \pi_\ell^{***}(a,b^{\prime\prime}),b^\prime\rangle.
$$
It follows that
$$
L_a^{**}b^{\prime\prime}=\pi_\ell^{**}(a,b^{\prime\prime}).
$$
Let $(b^\prime_\alpha)_\alpha\subseteq B^*$ such that $b^\prime_\alpha\stackrel{w^*} {\rightarrow}b^{\prime}$.
Since $\pi_\ell^{***}(a,b^{\prime\prime})\in B$, we have
$$
 \langle  b^{\prime\prime},L^*_a b_\alpha^\prime\rangle= \langle  L^{**}_ab^{\prime\prime}, b_\alpha^\prime\rangle= \langle  \pi_\ell^{***}(a,b^{\prime\prime}),b_\alpha^\prime\rangle=\langle  \pi_\ell^{***}(a,b^{\prime\prime}),b^\prime\rangle= \langle  b^{\prime\prime},L^*_a b^\prime\rangle.
$$
We conclude that $L_a^*$ is $weak^*-to-weak$ continuous, so $L_a$ is weakly compact.\\
Conversely, assume that $b^{\prime\prime}\in B^{**}$. Then by  $Goldstine^,s$ theorem [9, P.424-425], there is a net $(b_\alpha)_\alpha\subseteq B$ such that $b_\alpha\stackrel{w^*} {\rightarrow}b^{\prime\prime}$. Since for all $a\in A$, the operator $L_a$ is weakly compact, there is a subnet   $(b_{\alpha_\beta})_{_\beta}$ from $(b_\alpha)_\alpha$ such that $(L_a(b_{\alpha_\beta}))_{\beta}$ is weakly convergence to some point of $B$. Since $b_\alpha\stackrel{w^*} {\rightarrow}b^{\prime\prime}$,  $(L_a(b_{\alpha_\beta})=\pi(a,b_{\alpha_\beta}))_\beta$ is weakly convergence to $\pi_\ell^{***}(a,b^{\prime\prime})$. Consequently, for all $b^\prime\in B^*$, we have
$$
 \langle  \pi_\ell^{***}(a,b^{\prime\prime}),b^\prime\rangle=lim_\beta \langle  b^\prime,\pi(a,b_{\alpha_\beta})\rangle
=lim_\beta \langle  b^\prime,L_ab_{\alpha_\beta})\rangle.
$$
It follows that $\pi_\ell^{***}(a,b^{\prime\prime})\in B$.\\
\end{proof}

\noindent For a right Banach $A-module$ $B$, we can write Theorem 2-16 as follows.\\
Assume for all $a\in A$,  $R_a$ be the linear mapping from $B$ into itself such that $R_ab=\pi_r(b,a)$ for all $b\in B$. Then $B^{**}A\subseteq B$ if and only if $R_a$ is weakly compact.
Proof of this assertion is similar to proof of Theorem 2-16.\\\\

\noindent {\it{\bf Corollary 2-17.}}\\ i) Suppose that $B$ is a left Banach $A-module$. Then $AB^{**}\subseteq B$ if and only if $B$ is $Lwcc$.\\
ii) Suppose that $B$ is a right  Banach $A-module$. Then $B^{**}A\subseteq B$ if and only if $B$ is $Rwcc$.\\\\

\noindent {\it{\bf Corollary 2-18.}} Let $A$ be a $WSC$ Banach algebra with a $BAI$. If $A$ is Arens regular and  $A$ is a left ideal in its second dual, then $A$ is reflexive.
\begin{proof} Since $A$ is a left ideal in $A^{**}$, by using proceeding corollary, $A$ is $Lwcc$. Then by using Corollary 2.8 from [16], we are done.\\
\end{proof}

\noindent {\it{\bf Example 2-19.}}\\ i) Let $G$ be a compact group. Then we know that $L^1(G)$ is a left ideal $L^1(G)^{**}$ (resp. $M(G)^{**}$), and so by Corollary 2-18 (resp. Corollary 2-17), $L^1(G)$ (resp. $M(G)$) is a $Lwcc$. \\
ii) Corollary 2-18 shows that if $G$ is a finite group. Then $L^1(G)$ and $M(G)$ are reflexive.\\
iii) Let $G$ be a locally compact Hausdorff group. Let $X$ be a subsemigroup of $G$ which is the clouser of an open subset, and which contains the identity $e$ of $G$. Let $Z$ be a closed two-sided proper ideal in $X$ with the property that $X \setminus Z$ is relatively compact. Let $S$ be the quotient of $X$  obtained
by identifying all points of $Z$, more formally, for $x,y\in X$ write $x\sim y$ if either $x=y$ or both $x\in Z$ and $y\in Z$ and write $S=X/\sim$. By using Corollary 3.3 from [21],  $L^1(S)$ is an ideal $L^1(S)^{**}$, and so by using  Corollary 2-17,  $L^1(S)$ is a $Lwcc$.\\\\

\noindent {\it{\bf Theorem 2-20.}} Let $A$ be a $WSC$ Banach algebra with a $BAI$. If $A$ is Arens regular and  $A$ is a right ideal in its second dual, then $A$ is reflexive.
\begin{proof} Proof is similar to Corollary 2-18.\\ \end{proof}

\noindent {\it{\bf Definition 2-21.}} Suppose that $B$ is a left Banach $A-module$. Let $(e_{\alpha})_{\alpha}\subseteq A$ be  left approximate identity for $A$. We say that $(e_{\alpha})_{\alpha}$ is $weak^*$ left approximate identity ($=W^*LAI$) for $B^*$, if for all $b^\prime\in B^*$, we have $\pi_\ell (e_{\alpha},b^\prime) \stackrel{w^*} {\rightarrow}
b^\prime$. The definition of the $weak^*$ right approximate identity ($=W^*RAI$) is similar.\\
 We say that $(e_{\alpha})_{\alpha}$ is a $weak^*$ approximate identity  ($=W^*AI$) for $B^*$, if $B^*$ has $weak^*$ left and right approximate identity   that are equal.\\\\

\"{U}lger in [22] shows that for a Banach algebra $A$ with a $BAI$, if $A$ is a bisded ideal in its second dual, then $AA^*=A^*A$ and if $A$ is Arens regular, then $A^*$ factors on the both side. In the following, we extend these problems for module actions with some results in group algebras.\\
Let $B$ be a left Banach $A-module$. Then, $b^\prime\in B^*$ is said to be left weakly almost periodic functional if the set $\{\pi_\ell(b^\prime,a):~a\in A,~\parallel a\parallel\leq 1\}$ is relatively weakly compact. We denote by $wap_\ell(B)$ the closed subspace of $B^*$ consisting of all the left weakly almost periodic functionals in $B^*$.\\
The definition of the right weakly almost periodic functional ($=wap_r(B)$) is the same.\\
By [5, 16, 20], the definition of $wap_\ell(B)$ is equivalent to the following $$ \langle  \pi_\ell^{***}(a^{\prime\prime},b^{\prime\prime}),b^\prime\rangle=
 \langle  \pi_\ell^{t***t}(a^{\prime\prime},b^{\prime\prime}),b^\prime\rangle$$
for all $a^{\prime\prime}\in A^{**}$ and $b^{\prime\prime}\in B^{**}$.
Thus, we can write \\
$$wap_\ell(B)=\{ b^\prime\in B^*:~ \langle  \pi_\ell^{***}(a^{\prime\prime},b^{\prime\prime}),b^\prime\rangle=
 \langle  \pi_\ell^{t***t}(a^{\prime\prime},b^{\prime\prime}),b^\prime\rangle~~$$$$for~~all~~a^{\prime\prime}\in A^{**},~b^{\prime\prime}\in B^{**}\}.$$\\

By using [20], $b^\prime\in wap_\ell(B)$ if and only if for each sequence $(a_n)_n\subseteq A$ and $(b_m)_m\subseteq B$ and each $b^\prime\in B^*$, we have
    $$\lim_m\lim_n \langle  b^\prime,\pi_\ell(a_n,b_m)\rangle=\lim_n\lim_m \langle  b^\prime,\pi_\ell(a_n,b_m)\rangle,$$
    whenever both the iterated limits exist.\\
It is clear that $wap_\ell(B)=B^*$ if and only if ${Z}^\ell_{A^{**}}(B^{**})=B^{**}$.\\\\

\noindent{\it{\bf Definition 2-22.}} Let $B$ be a left Banach  $A-module$ and $A$ has a $BAI$ as $(e_\alpha)_\alpha$. We introduce the following subspace of $B^*$.
$$\ell(B^*)=\{b^\prime\in B^*:~\pi_\ell^*(b^\prime,e_\alpha)\stackrel{w} {\rightarrow} b^\prime\}.$$
Let $B$ be a right Banach  $A-module$ and $A$ has a $BAI$ as $(e_\alpha)_\alpha$.  Such as  proceeding definition, we introduce the following subspace of $B^*$.
$$\Re(B^*)=\{b^\prime\in B^*:~\pi_r^{t*}(b^\prime,e_\alpha)\stackrel{w} {\rightarrow} b^\prime\}.$$\\
If $\ell(B^*)=B^*$ (resp. $\Re(B^*)=B^*$), then it is clear that $(e_\alpha)_\alpha\subseteq A$ is a weakly right (resp. left) approximate identity for $B^{*}$. Therefore by using Lemma 2-8, $\ell(B^*)=B^*$ (resp. $\Re(B^*)=AB^*$) if and only if $B^*$ factors on the left (resp. right).\\\\

\noindent{\it{\bf Theorem 2-23.}} Let $B$ be a left Banach  $A-module$ and $A$ has a $RBAI$ as $(e_\alpha)_\alpha$. Then we have the following assertions.\\
i) $\ell(B^*)=B^*A$.\\
ii) $wap_\ell(B)\subseteq \ell(B^*)$, if $B^*$ has $W^*LAI$ as $A-module$ $(e_\alpha)_\alpha$.\\
\begin{proof} i) Let $\pi_\ell:A\times B\rightarrow B$ be the left module action such that $\pi_\ell(a,b)=ab$ for all $a\in A$ and $b\in B$. Thus for every $a\in A$, $b^\prime\in B^*$ and $b^{\prime\prime}\in B^{**}$, we have
$$ \langle  b^{\prime\prime},\pi_\ell^*(b^\prime a,e_\alpha)\rangle= \langle  b^{\prime\prime},\pi_\ell^*(b^\prime ,ae_\alpha)\rangle=
 \langle  \pi_\ell^{**}(b^{\prime\prime},b^\prime) ,ae_\alpha\rangle\rightarrow  \langle  \pi_\ell^{**}(b^{\prime\prime},b^\prime) ,a\rangle$$
$$= \langle  b^{\prime\prime},\pi_\ell^*(b^\prime ,a)\rangle= \langle  b^{\prime\prime},b^\prime a\rangle.$$
It follow that $\pi_\ell^*(b^\prime a,e_\alpha)\stackrel{w} {\rightarrow}b^\prime a$ and so $b^\prime a\in\ell(B^*)$.
 For  reverse inclusion, since by $Cohen^,s$ factorization theorem, we have $B^*A$ is a closed subspace of $B^*$, $\ell(B^*)\subseteq B^*A$.\\
 ii) Let $b^\prime\in wap_\ell(B)$. Since $(e_\alpha)_\alpha$ is $W^*LAI$ for $B^*$, $\pi_\ell^*(b^\prime,e_\alpha)\stackrel{w^*} {\rightarrow} b^\prime$. Also the set $\{\pi_\ell^*(b^\prime,e_\alpha):~\alpha\in I\}$ is relatively weakly compact which implies that
$\pi_\ell^*(b^\prime,e_\alpha)\stackrel{w} {\rightarrow} b^\prime$.\\
\end{proof}
\noindent{\it{\bf Corollary 2-24.}} Let $B$ be a left Banach  $A-module$ and $A$ has a $RBAI$ and $Z_{A^{**}}(B^{**})=B^{**}$. If $B^*$ has $W^*LAI$, then $B^*$ factors on the left.\\\\

\noindent{\it{\bf Example 2-25.}}\\ i) Let $G$ be a finite group. Then, by using proceeding corollary, we conclude that  $L^\infty(G)M(G)=RUC(G)=L^\infty(G)$.\\
ii) Let $G$ be an infinite compact group. Since $L^1(G)$ has a $BAI$, $L^\infty (G)$ has a $W^*BAI$. By using Proposition 4.4  from [22], we know that $wap(L^1(G))=C(G)$. Thus, by using proceeding theorem, $C(G)\subseteq \ell (L^\infty (G))$.\\\\

\noindent{\it{\bf Theorem 2-26.}} Let $B$ be a right Banach  $A-module$ and $A$ has a $LBAI$ as $(e_\alpha)_\alpha$. Then we have the following assertions.\\
i) $\Re(B^*)=AB^*$.\\
ii) $wap_r(B)\subseteq \Re(B^*)$, if $B^*$ has $W^*RAI$ $A-module$ $(e_\alpha)_\alpha$.
\begin{proof} Proof is similar to Theorem 2-23.\\
\end{proof}

\noindent{\it{\bf Corollary 2-27.}} Let $B$ be a right Banach  $A-module$ and $A$ has a $BAI$ as $(e_\alpha)_\alpha$.
If $B$ factors on the left, then $wap_r(B)\subseteq\Re (B^*)$.\\
\begin{proof} Let $b\in B$ and $b^\prime \in B^*$. Since $B$ factors on the left, there are $a\in A$ and $y\in B$ such that $b=ya$. Then
$$ \langle  \pi_r^{t*}(b^\prime,e_\alpha),b\rangle= \langle  b^\prime,\pi_r^{t}(e_\alpha,b)\rangle= \langle  b^\prime,\pi_r(b,e_\alpha)\rangle=
 \langle  b^\prime,\pi_r(ya,e_\alpha)\rangle$$$$= \langle  \pi_r^*(b^\prime,y),ae_\alpha\rangle\rightarrow
 \langle  \pi_r^*(b^\prime,y),a\rangle= \langle  b^\prime,ya\rangle= \langle  b^\prime, b\rangle.$$
It follows that $\pi_r^{t*}(b^\prime,e_\alpha)\stackrel{w^*} {\rightarrow}b^\prime $. Then by using Theorem 2-26, we are done. \\
\end{proof}

\noindent{\it{\bf Corollary 2-28.}} Let $B$ be a right Banach  $A-module$ and $A$ has a $BAI$ and $Z^t_{B^{**}}(A^{**})=B^{**}$. If $B^*$ has $W^*RAI$ $A-module$, then $B^*$ factors on the right.\\\\

\noindent{\it{\bf Theorem 2-29.}} We have the following statements.
\begin{enumerate}
\item Let $B$ be a left Banach  $A-module$. If $AB^{**}\subseteq B$, then $B^*A\subseteq wap_\ell(B)$.
\item Let $B$ be a right Banach  $A-module$. If $B^{**}A\subseteq B$, then $AB^*\subseteq wap_r(B)$.
\end{enumerate}
\begin{proof}
\begin{enumerate}
\item By Theorem 2-17, we know that $B$ is $Lwcc$. Let $a\in A$ and suppose that $L_a$ is the mapping from $B$ into itself by definition $L_a(b)=\pi_\ell(a,b)$ for each $b\in B$. By easy calculation, it is clear that $(L_a)^*(b^\prime)=\pi^*_\ell(b^\prime,a)$. Since $L_a$ is weakly compact, $(L_a)^*$ is weakly compact.
Then the set $$\{(L_a)^*(\pi^*_\ell(b^\prime,x)):~x\in A_1\},$$ is weakly compact. Now let $x\in A_1$ and $y\in B$. Then we have the following equality
$$ \langle  (L_a)^*(\pi^*_\ell(b^\prime,x)),y\rangle= \langle  \pi^*_\ell(b^\prime,x),L_a(y)\rangle= \langle  \pi^*_\ell(b^\prime,x),\pi_\ell(a,y)\rangle$$
$$= \langle  \pi_\ell^*(\pi^*_\ell(b^\prime,x),a),y)\rangle.$$
 It follows that the mapping $\pi_\ell^*(\pi^*_\ell(b^\prime,x),a)$ is weakly compact for each $a\in A$ and $b^\prime\in  B^*$. Hence $\pi_\ell^*(b^\prime,x)\in wap_\ell(B)$, and so $B^*A\subseteq wap_\ell(B)$.
\item Proof is similar to  proceeding proof.\\
\end{enumerate}
\end{proof}

\noindent {\it{\bf Example 2-30.}}\\  i) Let $G$ be a locally compact group and $1\leq p\leq \infty$. We know that $L^p(G)$ is the left Banach $L^1(G)-module$ under  convolution as multiplication. Assume that for all $f\in L^1(G)$, $L_f:L^p(G)\rightarrow L^p(G)$ be the linear mapping such that $L_fg=f*g$ whenever $g\in L^p(G)$. Then, since $L^1(G)L^p(G)^{**}= L^1(G)L^p(G)\subseteq  L^p(G)$ for all $1 <  p <  \infty$, by Theorem 2-28, $L_f$ is weakly compact.\\
It is the same that for all $\mu\in M(G)$, the mapping $L_\mu$ from $L^p(G)$ into itself with $L_\mu f=\mu*f$ is weakly compact whenever $1 <   p <  \infty$.\\
ii) Let $G$ be an infinite compact group. Then we know that $L^1(G)$ is an ideal in its second dual, $L^1(G)^{**}$. Therefore, by using proceeding theorem and Proposition 3.3, from [22], we have $LUC(G)=L^\infty (G)L^1(G)\subseteq wap(L^1(G))$ and $RUC(G)=L^1(G)L^\infty(G)\subseteq wap(L^1(G))$. By using Proposition 4.4 from [22], since $wap(L^1(G))=C(G)$, we conclude that $LUC(G)\cap RUC(G)\subseteq C(G)$, and so $LUC(G)\cap RUC(G)=C(G)$.\\\\

Let $B$ be a Banach  $A-bimodule$ and $a^{\prime\prime}\in A^{**}$. We define the locally topological centers of the left and right module actions of $a^{\prime\prime}$ on $B$, respectively, as follows\\
$$Z_{a^{\prime\prime}}^t(B^{**})=Z_{a^{\prime\prime}}^t(\pi_\ell^t)=\{b^{\prime\prime}\in B^{**}:~~~\pi^{t***t}_\ell(a^{\prime\prime},b^{\prime\prime})=
\pi^{***}_\ell(a^{\prime\prime},b^{\prime\prime})\},$$
$$Z_{a^{\prime\prime}}(B^{**})=Z_{a^{\prime\prime}}(\pi_r^t)=\{b^{\prime\prime}\in B^{**}:~~~\pi^{t***t}_r(b^{\prime\prime},a^{\prime\prime})=
\pi^{***}_r(b^{\prime\prime},a^{\prime\prime})\}.$$\\
It is clear that ~~~~~~~$$\bigcap_{a^{\prime\prime}\in A^{**}}Z_{a^{\prime\prime}}^t(B^{**})=Z_{A^{**}}^t(B^{**})=
Z(\pi_\ell^t),$$ ~~~~~~~ ~~\\ $$~~~~~~~~~~~~~~~~~~~~~~~~~~~~\bigcap_{a^{\prime\prime}\in A^{**}}Z_{a^{\prime\prime}}(B^{**})=Z_{A^{**}}(B^{**})=
Z(\pi_r).$$\\
The definition of $Z_{b^{\prime\prime}}^t(A^{**})$ and $Z_{b^{\prime\prime}}(A^{**})$ for some $b^{\prime\prime}\in B^{**}$ are the same.\\\\

\noindent{\it{\bf Theorem 2-31.}} Let $B$ be a Banach left $A-module$ and $A$ has a $LBAI$ $(e_{\alpha})_{\alpha}\subseteq A$ such that $e_{\alpha} \stackrel{w^*} {\rightarrow}e^{\prime\prime}$ in $A^{**}$ where $e^{\prime\prime}$ is a left unit for $A^{**}$. Suppose that $Z^t_{e^{\prime\prime}}(B^{**})=B^{**}$. {Then, $B$ factors on the right with respect to $A$ if and only if
$e^{\prime\prime}$ is a left unit for $B^{**}$.
\begin{proof}
Assume that  $B$ factors on the right with respect to $A$. Then  for every $b\in B$, there are $x\in B$ and $a\in A$ such that $b= ax$. Then for every $b^\prime\in B^*$, we have
$$\langle  \pi^*_\ell(b^\prime,e_\alpha),b\rangle =\langle  b^\prime,\pi_\ell(e_\alpha,b)\rangle =\langle  \pi^{***}_\ell(e_\alpha,b),b^\prime\rangle $$$$=
\langle  \pi^{***}_\ell(e_\alpha,ax),b^\prime\rangle =
\langle  \pi^{***}_\ell(e_\alpha a,x),b^\prime\rangle $$$$=\langle  e_\alpha a,\pi^{**}_\ell(x,b^\prime)\rangle =
\langle  \pi^{**}_\ell(x,b^\prime),e_\alpha a\rangle $$$$\rightarrow \langle   \pi^{**}_\ell(x,b^\prime),a\rangle =\langle  b^\prime,b\rangle .$$
It follows that  $\pi^*_\ell(b^\prime,e_\alpha) \stackrel{w^*} {\rightarrow}b^\prime$ in $B^*$.
Let $b^{\prime\prime}\in B^{**}$ and $(b_\beta)_\beta\subseteq B$ such that $b_\beta\stackrel{w^*} {\rightarrow}b^{\prime\prime}$ in $B^{**}$. Since $Z^t_{e^{\prime\prime}}(B^{**})=B^{**}$, for every $b^\prime\in B^*$, we have the following equality
$$\langle  \pi_\ell^{***}(e^{\prime\prime},b^{\prime\prime}),b^\prime\rangle =
\lim_\alpha\lim_\beta\langle  b^\prime, \pi_\ell(e_\alpha,b_\beta)\rangle $$$$=
\lim_\beta\lim_\alpha\langle  b^\prime, \pi_\ell(e_\alpha,b_\beta)\rangle =
\lim_\beta\langle  b^\prime,b_\beta\rangle $$$$
=\langle  b^{\prime\prime},b^\prime\rangle .$$
It follows that $\pi_\ell^{***}(e^{\prime\prime},b^{\prime\prime})=b^{\prime\prime}$, and so $e^{\prime\prime}$ is a left unit for $B^{**}$.\\
Conversely, let $e^{\prime\prime}$ be a left unit for $B^{**}$ and suppose that $b\in B$. Thren for every $b^\prime\in B^*$, we have
$$\langle  b^\prime,\pi(e_\alpha,b)\rangle =\langle  \pi^{***}(e_\alpha,b),b^\prime\rangle =\langle  e_\alpha,\pi^{**}(b,b^\prime)\rangle
=\langle  \pi^{**}(b,b^\prime),e_\alpha\rangle $$$$=
\langle  e^{\prime\prime},\pi^{**}(b,b^\prime)\rangle =\langle  \pi^{***}(e^{\prime\prime},b),b^\prime\rangle =\langle  b^\prime,b\rangle .$$
Then we have $\pi^*_\ell(b^\prime,e_\alpha) \stackrel{w} {\rightarrow}b^\prime$ in $B^*$, and so by Cohen factorization theorem we are done.

\end{proof}

\noindent{\it{\bf Corollary 2-32.}} Let $B$ be a Banach left $A-module$ and $A$ has a $LBAI$ $(e_{\alpha})_{\alpha}\subseteq A$ such that $e_{\alpha} \stackrel{w^*} {\rightarrow}e^{\prime\prime}$ in $A^{**}$ where $e^{\prime\prime}$ is a left unit for $A^{**}$. Suppose that $Z^t_{e^{\prime\prime}}(B^{**})=B^{**}$. Then
$\pi^*_\ell(b^\prime,e_\alpha) \stackrel{w} {\rightarrow}b^\prime$ in $B^*$ if and only if $e^{\prime\prime}$ is a left unit for $B^{**}$.\\\\

Hu, Neufang and Ruan in [15], have been studied multiplier on new class of Banach algebras. They showed that how a multiplier on Banach algebra $A$ to be implemented by an element from $A$ is determined by its behavior on $A^*$ and $A^{**}$, respectively. In the following we study some of these discussion on module actions with some results, see [15, Theorem 3].\\
 For a Banach algebra $A$, we recall that a bounded linear operator $T:A\rightarrow A$ is said to be a left (resp. right) multiplier
if, for all $a, b\in A$, $T(ab)=T(a)b$ (resp. $T(ab)=aT(b)$). We denote by $LM(A)$ (resp. $RM(A)$) the set of all left (resp. right) multipliers of $A$. The set $LM(A)$ (resp. $RM(A)$) is normed subalgebra of the algebra $L(A)$ of bounded linear operator on $A$.\\
 Let  $B$ be a Banach left [resp. right] $A-module$ and  $T\in \mathbf{B}(A,B)$. Then $T$ is called extended left [resp. right] multiplier if ~$T(a_1 a_2)=\pi_r(T(a_1),a_2)\\ ~[resp.~T(a_1 a_2)=\pi_\ell(a_1,T(a_2))]$ ~for all $a_1,a_2\in A$.\\
We show by $LM(A,B)~[resp. ~ RM(A,B)]$ all of the Left [resp. right] multiplier extension from $A$ into $B$. \\
 Let $a^\prime\in A^*$. Then the mapping $T_{a^\prime}:a\rightarrow {a^\prime}a ~[resp.~~ R_{a^\prime}~a\rightarrow a{a^\prime}]$ from $A$ into $A^*$ is left [right] multiplier, that is, $T_{a^\prime}\in LM(A,A^*)~[R_{a^\prime}\in RM(A,A^*)].$ $T_{a^\prime}$ is weakly compact if and only if ${a^\prime}\in wap(A)$. So, we can write $wap(A)$ as a subspace of $LM(A,A^*)$.\\\\

\noindent{\it{\bf Theorem 2-33.}} Let $B$ be a   Banach $A-bimodule$ with a $BAI$ $(e_\alpha)_\alpha \subseteq A$. Then \begin{enumerate}
\item If $T\in LM(A,B)$, then $T(a)=\pi_r^{***}(b^{\prime\prime},a)$ for some $b^{\prime\prime}\in B^{**}$.
\item If $T\in RM(A,B)$, then $T(a)=\pi_\ell^{***}(a,b^{\prime\prime})$ for some $b^{\prime\prime}\in B^{**}$.
\end{enumerate}
\begin{proof}
\begin{enumerate}
 \item Since $(T(e_\alpha))_\alpha \subseteq B$ is bounded, it has weakly limit point in $B^{**}$. Let $b^{\prime\prime}\in B^{**}$ be a weakly limit point of $(T(e_\alpha))_\alpha$ and without loss generally, take  $T(e_\alpha) \stackrel{w} {\rightarrow}b^{\prime\prime}$. Then for every $b^\prime\in B^*$ and $a\in A$, we have

$$\langle  \pi_r^{***}(b^{\prime\prime},a),b^\prime\rangle
=\lim_\alpha\langle  b^\prime, T(e_\alpha)a\rangle =\lim_\alpha\langle  b^\prime, T(e_\alpha a)\rangle $$$$=\lim_\alpha\langle  T^*(b^\prime), e_\alpha a\rangle =
\langle  T^*(b^\prime),  a\rangle =\langle  b^\prime,  T(a)\rangle .$$
It follows that $\pi_r^{***}(b^{\prime\prime},a)=T(a)$.
\item  Proof is similar to (1).\\
\end{enumerate}\end{proof}

In the proceeding theorem, if we take $B=A$, then we have the following statements
\begin{enumerate}
\item If $T\in LM(A)$, then $T(a)=a^{\prime\prime}a$  for some $a^{\prime\prime}\in A^{**}$.
\item If $T\in RM(A)$, then $T(a)=aa^{\prime\prime}$  for some $a^{\prime\prime}\in A^{**}$.\\
\end{enumerate}

\noindent{\it{\bf Definition 2-34.}} Let $B$ be a Banach left  $A-module$ and $b^{\prime\prime}\in B^{**}$. Suppose that $(b_{\alpha})_{\alpha}\subseteq B$ such that  $b_{\alpha} \stackrel{w^*} {\rightarrow}b^{\prime\prime}$. We define the following set
$$\widetilde{Z}_{b^{\prime\prime}}(A^{**})=\{a^{\prime\prime}\in A^{**}:~~~\pi_\ell^{***}(a^{\prime\prime},b_{\alpha})\stackrel{w^*} {\rightarrow}\pi_\ell^{***}(a^{\prime\prime},b^{\prime\prime})\},$$
which is subspace of $A^{**}$. It is clear that $Z_{b^{\prime\prime}}(A^{**})\subseteq \widetilde{Z}_{b^{\prime\prime}}(A^{**})$, and so
$$Z_{B^{**}}(A^{**})=\bigcap_{b^{\prime\prime}\in B^{**}}{Z}_{b^{\prime\prime}}(A^{**})\subseteq \bigcap_{b^{\prime\prime}\in B^{**}}\widetilde{Z}_{b^{\prime\prime}}(A^{**}).$$
For a Banach right $A-module$, the definition of $\widetilde{Z}^t_{a^{\prime\prime}}(B^{**})$  is similar.\\\\

\noindent{\it{\bf Theorem 2-35.}} Let $B$ be a left  Banach $A-module$ and $T\in \mathbf{B}(A,B)$. Consider the following statements.
\begin{enumerate}
\item $T=\ell_b$, for some $b\in B$.
\item  $T^{**}(a^{\prime\prime})=\pi_\ell^{***}(a^{\prime\prime},b^{\prime\prime})$ for some $b^{\prime\prime}\in B^{**}$ such that $\widetilde{Z}_{b^{\prime\prime}}(A^{**})=A^{**}$.
\item $T^*(B^*)\subseteq BB^*$.
\end{enumerate}
Then $(1)\Rightarrow (2)\Rightarrow (3)$.\\
Assume that $B$ has $WSC$. If we take $T\in RM(A,B)$ and  $B$ has a sequential $BAI$, then (1), (2) and (3) are equivalent.\\
\begin{proof}
$(1)\Rightarrow (2)$\\
 Let $T=\ell_b$, for some $b\in B$. Then $T^{**}(a^{\prime\prime})=\ell^{**}_b(a^{\prime\prime})=\pi_\ell^{***}(a^{\prime\prime},b)$ for every $a^{\prime\prime}\in A^{**}$, and so proof is hold.\\
 $(2)\Leftrightarrow (3)$\\
Take $a^{\prime\prime}\in (BB^*)^\bot$. Assume that $b^{\prime\prime}\in B^{**}$ and $(b_{\alpha})_{\alpha}\subseteq B$ such that  $b_{\alpha} \stackrel{w^*} {\rightarrow}b^{\prime\prime}$. For every $b^\prime\in B^{**}$, we have the following equality
$$\langle  a^{\prime\prime}, T^*(b^\prime)\rangle =\langle  T^{**}(a^{\prime\prime}), b^\prime\rangle =\langle  \pi_\ell^{***}(a^{\prime\prime},b^{\prime\prime}), b^\prime\rangle =\lim_\alpha\langle  \pi_\ell^{***}(a^{\prime\prime},b_{\alpha}),b^{\prime}\rangle $$$$=
\lim_\alpha\langle  a^{\prime\prime},\pi_\ell^{**}(b_{\alpha},b^{\prime})\rangle =0.$$
It follows that $T^*(B^*)\subseteq BB^*$.\\
Take $T\in RM(A,B)$ and suppose that  $B$ is $WSC$ with sequential $BAI$. It is enough, we show that $(3)\Rightarrow (1)$.
Assume that $(e_n)_n\subseteq A$ is a $BAI$ for $B$. Then for every $b^\prime\in B^*$, we have
$$\mid\langle  b^\prime,T(e_n)\rangle -\langle  b^\prime,T(e_m)\rangle \mid=\mid \langle  T^*(b^\prime),e_n-e_m\rangle \mid=
\mid \langle  \pi^{**}_\ell(b,b^\prime),e_n-e_m\rangle \mid$$$$=
\mid \langle  b,\pi^{*}_\ell(b^\prime,e_n-e_m)\rangle \mid=
\mid \langle  b^\prime,\pi^{}_\ell(e_n-e_m,b)\rangle \mid\rightarrow 0.$$
It follows that $(T(e_n))_n$ is weakly Cauchy sequence in $B$ and since $B$ is $WSC$, there is $b\in B$ such that
$T(e_n)\stackrel{w} {\rightarrow}b$ in $B$. Let $a\in A$. Then for every $b^\prime\in B^*$, we have
$$\langle  b^\prime ,\pi_\ell(a,b)\rangle =\langle  \pi_\ell^*(b^\prime ,a),b)\rangle =\lim_n\langle  \pi_\ell^*(b^\prime ,a),T(e_n)\rangle $$$$
=\lim_n\langle  b^\prime ,\pi_\ell(a,T(e_n))\rangle =\lim_n\langle  b^\prime ,T(ae_n)\rangle $$$$=\lim_n\langle  T^*(b^\prime),ae_n\rangle
=\langle  T^*(b^\prime),a\rangle $$$$=\langle  b^\prime,T(a)\rangle .$$
Thus $\ell_b(a)=\pi_\ell(a,b)=T(a)$.\\
\end{proof}

\noindent{\it{\bf Example 2-36.}} Let $G$ be a locally compact group. Then by convolution multiplication,  $M(G)$ is a  $L^1(G)-bimodule$. Let $f\in L^1(G)$ and $T(\mu)=\mu*f$ for all $\mu\in M(G)$. Then $T^*(L^\infty (G))\subseteq M(G)M(G)^*$. Also if we take  $T(\mu)=f*\mu$ for all $\mu\in M(G)$, then we have $T^*(L^\infty (G))\subseteq M(G)^*M(G)$.\\

\bibliographystyle{amsplain}

\it{Department of Mathematics, University of Mohghegh Ardabili, Ardabil, Iran\\
{\it Email address:} haghnejad@aut.ac.ir\\\\

\end{document}